\documentclass[reviewcopy]{elsart}

\usepackage{latexsym}
\usepackage{amsmath}
\usepackage{amscd}
\usepackage[all]{xy} %This is XY-pic.
\usepackage{amsfonts}
\usepackage{amssymb}
\usepackage{amsbsy}

\journal{\mbox{}}

\newcommand{\ot}{\otimes}
\newcommand{\ra}{\rightarrow}

\newcommand{\BC}{\mathbb{C}}
\newcommand{\BR}{\mathbb{R}}

\newcommand{\BZ}{\mathbb{Z}}

\newcommand{\dirac}{D\!\!\!\!/\,}

\newcommand{\I}{\mathbf{i}}

\newcommand{\Spinc}{\ensuremath{\mathrm{Spin}^c}}
\newcommand{\Kahler}{K\"{a}hler}

\newcommand{\T}{{\ensuremath{\textsc{T}}}}

\newtheorem{theorem}{Theorem}

\newif\ifpdf
\ifx\pdfoutput\undefined
\pdffalse % we are not running PDFLaTeX
\else
\pdfoutput=1 % we are running PDFLaTeX
\pdftrue \fi

\ifpdf
\usepackage[pdftex]{graphicx}
\usepackage{hyperref}
\else %\topmargin .5in
\usepackage{graphics}
\usepackage{psfrag}
\fi

\begin{document}

\ifpdf \DeclareGraphicsExtensions{.jpg,.pdf,.mps,.png} \else
\DeclareGraphicsExtensions{.eps} \fi

%Declares which type of graphics extension based upon whether this is a pdf run
%or a .ps run.

\begin{frontmatter}


%\author{Scott Jeremy Baldridge\corauthref{cor1}\thanksref{label2}}


\title{New symplectic $4$--manifolds with $b_+{=}1$}

\author{Scott Baldridge}%\corauthref{cor1}\thanksref{label2}}
\address{Department of Mathematics,  Indiana University,
 Bloomington, Indiana 47405, USA,  FAX:(812) 855-0046 }

%\ead{sbaldrid@indiana.edu}

%**************************************************
%******* Abstract     *****************************
%**************************************************

\begin{abstract}
Symplectic 4-manifolds $(X,\omega)$ with $b_+{=}1$ are roughly
classified by the canonical class $K$ and the symplectic form
$\omega$ depending upon the sign of $K^2$ and $K\cdot \omega$.
Examples are known for each category except for the case when the
manifold satisfies $K^2=0$, $K\cdot \omega >0$, $b_1=2$, and fails
to be of Lefschetz type.  The purpose of this paper is to
construct an infinite number of examples of such manifolds.
Furthermore, we will show that these manifolds have very special
properties --- they are not complex manifolds, their
Seiberg-Witten invariants are independent of the chamber
structure, and  they do not have metrics of positive scalar
curvature.

\end{abstract}


\end{frontmatter}

%**************************************************
%******* Introduction *****************************
%**************************************************

\section{Introduction}

Symplectic manifolds with $b_+=1$ are useful because they often
form the basic building blocks for symplectic manifolds with
$b_+>1$. Furthermore, the Seiberg-Witten invariants have a richer
structure when $b_+=1$ allowing one to make and test hypotheses
about symplectic manifolds that might be true generally.  In this
paper we describe new examples of symplectic manifolds with
$b_+=1$.

We begin the discussion with a quick treatment of the
classification of symplectic 4--manifolds and mention the history
behind the examples described in this paper. Symplectic
4--manifolds have a natural, well-defined extension of the Kodaira
dimension used to classify compact complex surfaces defined by the
sign of two numbers:  $$K^2 \hspace{.5cm} \mbox{ and }
\hspace{.5cm} K^2\cdot [\omega],$$ where $K$ is the canonical
class of a symplectic 4--manifold $(X,\omega)$
\cite{symp:li,symp:survey_symp_b_+=1}. In particular, for a
minimal symplectic 4--manifold,

$$\kappa(X,\omega) = \left\{ \begin{array}{cl} \ \  -\infty \ \  & \mbox{if $K^2<0$ or $K\cdot
[\omega]<0$}\\
0 & \mbox{if $K^2=0$, \ $K\cdot [\omega] =0$}\\
1 & \mbox{if $K^2=0$, \ $K\cdot [\omega] >0$}\\
2 & \mbox{if $K^2>0$ and $K\cdot [\omega] >0$.}\end{array}\right.
$$

For a general symplectic 4--manifold, the Kodaira dimension is
defined to be the Kodaira dimension of a minimal model of
$(X,\omega)$.  (The minimal model may not be unique, but
$\kappa(X,\omega)$ is still well-defined.) Note that examples with
$K^2>0$, $K\cdot \omega=0$ do not exist
\cite{symp:li,sw:sw_inv_and_symp_form}.

\Kahler\ surfaces provide many of the known examples of minimal
symplectic 4--manifolds.  For instance,  $\kappa(X,\omega) =
-\infty$ if and only if $X$ is diffeomorphic to a rational or
ruled surface \cite{symp:app_gen_wall_cross}.  Minimal symplectic
4--manifolds with $\kappa = 0$ include $K3$, Enriques surface, and
hyperelliptic surfaces. Incidently, the first example  of a
non-\Kahler\ symplectic 4--manifold satisfied $\kappa =0$
\cite{symp:symp_example}. Products $T^2\times \Sigma_g$ where
$g>1$ satisfy $\kappa(T^2\times \Sigma_g, \omega) =1$ and \Kahler\
surfaces of general type are examples of symplectic 4--manifolds
with $\kappa = 2$.

In this paper we are interested in symplectic 4--manifolds that
satisfy ${\kappa=1}$ and ${b_+=1}$. For these manifolds the first
Betti number is always zero  or two:  Noether's formula requires
that $b_1$ is even and the restriction on the size of $b_1$
follows from the Hirzebruch signature theorem,
$$0 = K^2 = 2 \chi + 3
\sigma = 9-4b_1 - b_{-},$$ where $\chi$ is the Euler class of $X$
and $\sigma$ is its signature.

The $b_1=0$ case is covered by Dolgachev surfaces. Recently, using
a construction of Fintushel and Stern, Park produced similar
non-complex symplectic $4$--manifolds (see
\cite{sw:non_cx_symp_4_man}).  The $b_1=2$ case is interesting
because there are no examples from \Kahler\ surfaces. McDuff and
Salamon mentioned this fact in a 1995 survey paper and went on to
discuss how one might construct symplectic manifolds for this case
\cite{symp:survey_symp_b_+=1}. They broke the search for such
manifolds into two cases based upon cup product structure of
$H^1(X;\BZ)$.  These two cases can be reformulated in terms of
\Kahler-like condition called Lefschetz type.

Symplectic 4--manifolds $(X,\omega)$ are said to be of {\em
Lefschetz type} if  ${[\omega]\in H^2(X;\BR)}$ satisfies the
conclusion of the Hard Lefschetz Theorem, namely, that $\cup
[\omega]:H^1(X;\BR) \ra H^3(X;\BR)$  is an isomorphism.
Essentially, McDuff and Salamon broke the search into 4--manifolds
which either did or did not have this condition.

Examples of $\kappa=1$, $b_+=1$, $b_1=2$ manifolds of Lefschetz
type can be constructed by first doing zero-surgery on a fibered
knot in $S^3$ to get a closed 3--manifold, then taking the product
with $S^1$.  The resulting 4--manifold then has the desired
properties. However, McDuff and Salamon pointed out in their
survey paper that the discovery of manifolds not of Lefschetz type
remained open.

At about the same time, Li and Liu discovered a general
wall-crossing formula for the Seiberg-Witten invariant when
$b_+{=}1$ \cite{sw:wallcross}. (Recall that the Seiberg-Witten
invariant generally depends upon the chamber in which it is
calculated when $b_+{=}1$.) In their paper Li and Liu questioned
whether there existed symplectic 4--manifolds where the
wall-crossing number was always equal to zero.  They were
interested in such examples because the Seiberg-Witten invariant
would depend only on the smooth structure, implying for instance
that these manifolds did not have any metrics with positive scalar
curvature. The only known examples were limited to special
$T^2$--bundles over $T^2$. These manifolds have Kodaira dimension
$\kappa=0$ (c.f. \cite{symp:t2bdls}) and one basic class (given by
the canonical \Spinc\ structure).

There is in fact a robust number of examples which satisfy the
conditions prescribed by both McDuff-Salamon and Li-Liu.  This
paper describes the conditions needed to produce them and provides
an infinite family of such examples.  In particular, we prove:

\begin{theorem} For each integer $g>1$, there exists
a smooth, closed, ${\kappa=1}$, $b_1{=}2$, $b_+{=}1$ symplectic
4--manifold $(X_g,\omega)$ not of Lefschetz type constructed in
\S2 below with the following properties:

\begin{enumerate}
\item $X_g$ does not admit a complex structure,

\item $X_g$ does not carry a metric with positive scalar
curvature, and

\item The Seiberg-Witten invariant is  a smooth invariant of $X_g$
given by
$$\mathcal{SW}_{X_g}(\mathfrak{s}) = (\mathfrak{s}^{-2}-3+\mathfrak{s}^2)^{g-1}.$$
In particular, $X_g$ is SW-simple type.
\end{enumerate}\label{thm:1}
\end{theorem}

When $g=1$, the manifold $X_g$ satisfies all of the conditions in
Theorem \ref{thm:1} except that $\kappa=0$.

\section{Construction}

In this section we construct a manifold $X_g$ with the properties
above for each integer $g>1$ using a technique similar to those
found in \cite{symp:symp_example,symp:sympmanfreecircleaction}.
Let $\Sigma$ be a surface of genus $g>0$.   Let $\langle
a_1,b_1,\ldots a_g,b_g\rangle$ be smooth loops which represent a
symplectic basis of $H_1(\Sigma;\BZ)$,
\begin{figure}[h]
\begin{center}
\psfrag{a1}{$a_1$} \psfrag{b1}{$b_1$}\psfrag{a2}{$a_2$}
\psfrag{a3}{$a_3$}\psfrag{b2}{$b_2$}\psfrag{b3}{$b_3$}
    \includegraphics{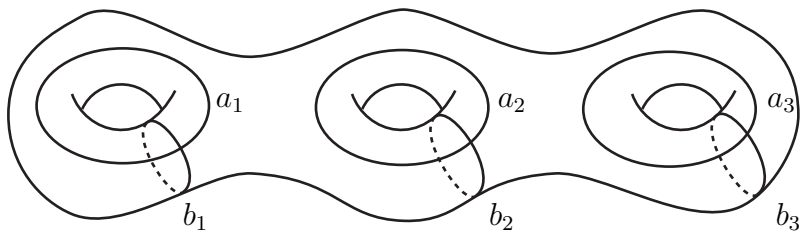}
    \caption{Surface of genus $3$. }
    \label{fig:surface_genus_g}
\end{center}
\end{figure}
and let $\langle \alpha_1, \beta_1, \ldots \alpha_g,
\beta_g\rangle$ be closed 1--forms which represent a dual basis
for $H^1(\Sigma;\BZ)$ with respect to the $a_i$'s and $b_i$'s,
i.e.,
$$\alpha_i(a_j) = \delta_{ij}, \hspace{1cm} \beta_i(b_j) = \delta_{ij}, \hspace{1cm} \mbox{ and
zero otherwise.}$$ Consider an orientation-preserving
diffeomorphism $\varphi:\Sigma \ra \Sigma$ with the following
matrix with respect to the basis above:

$$\varphi^* = \left(\begin{array}{cc} \begin{array}{cc} 1 & 0\\  1 & 1
\end{array}  & \hspace{.1cm} \mbox{\large 0} \hspace{.1cm}\\
\vspace{.2cm} \mbox{\large 0} & \mbox{\large A}
\end{array}\right)$$

where $\ker (A-Id) =0.$  Note that this means $\varphi^*[\beta_1]
= [\beta_1]$ and $\varphi^*[\alpha_1] = [\beta_1 + \alpha_1]$, and
that $\dim \ker (\varphi^* - Id) =1$.

Such diffeomorphisms exist for any genus.  For the purposes of
this paper, we will consider  diffeomorphisms given by the
following sequence of Dehn twists acting on the left,
$$
\varphi = (T_{b_g}T^{-1}_{a_g}) \cdots (T_{b_2}T^{-1}_{a_2}) \cdot
T_{a_1}
$$
for each genus $g$.
\bigskip

Next, create the mapping torus $Y$ obtained by crossing $\Sigma$
with the interval $[0,1]$ and identifying the ends by $\varphi$.
Then
$$Y = \left(\Sigma \times [0,1]\right) / \left((x,1) \sim (\varphi(x),0)\right)$$
is a smooth, closed, 3--dimensional manifold.

The 4--manifold that we will construct is a nontrivial  circle
bundle over $Y$.   This circle bundle has to be chosen carefully
because 4--manifolds which are circle bundles over a 3--manifold
which fibers over $S^1$ are not necessarily symplectic (cf.
\cite{sw:circleactions}).

In order to show that the 4--manifold has the desired properties
we need to describe its cohomology and the cup products
explicitly. We begin by first writing down a basis for the
cohomology of $Y$.

\bigskip

The  mapping torus $Y$ comes with a projection map, $p:Y\ra S^1$,
with fiber $\Sigma$. Let $dt$ be the volume form on $S^1$. Then
$\theta = p^*(dt)$ is a nowhere-zero closed 1--form in
$\Omega^1(Y;\BR)$. The form $\theta$ is one of the components
needed to construct the symplectic form on the 4--manifold.

There is another important 1--form on $Y$ to consider.  The
cohomology class $[\beta_1]$ is invariant under $\varphi^*$, so
there exists a function $f\in \Omega^0(\Sigma)$ such that
$\varphi^*(\beta_1) = \beta_1 + df$ point-wise. We can use this
1--form to construct a closed 1--form on $Y$ which represents a
nontrivial integral cohomology class. Let $\rho:[0,1] \ra [0,1]$
be a smooth function which is identically 0 near 0 and identically
1 near 1.  Extend $\beta_1$ to $\Sigma \times [0,1]$ by writing
$$\beta(x,t) = \rho(t)\beta_1(x) + (1-\rho(t))(\varphi^*(\beta_1(x)))
- (\frac{d}{dt}\rho(t))f(x)dt.$$ This glues up to give a smooth,
closed, 1--form on $Y$ since
$$\frac{d^n}{dt^n}\varphi^*(\beta)|_{\Sigma\times\{0\}} =
\frac{d^n}{dt^n}\beta|_{\Sigma\times\{1\}}$$ for all nonnegative
integers $n$ (we will also call the wrapped up 1-form $\beta$).

\bigskip

Next we show that the cohomology classes $[\theta]$ and $[\beta]$
form a basis for $H^1(Y;\BZ)$. It is clear that they are integral
and independent. Therefore we need only show that the dimension of
$H^1(Y)$ is 2, which follows from the Wang sequence and the fact
that $\dim \ker (\varphi^* -1) =1$.

\begin{eqnarray*}
\xymatrix{ H^0(\Sigma) \ar[r] \ar@{=}[d] & H^1(Y) \ar@{=}[d]
\ar[r]^{\ |_{\Sigma}} &H^1(\Sigma) \ar@{=}[d] \ar[r]^{\varphi^*-1}
& H^1(\Sigma)
\ar[r] \ar@{=}[d] & H^2(Y) \ar@{=}[d]\\
\BZ \ar[r] & \BZ\oplus\BZ  \ar[r] & \BZ \oplus \BZ^{2g-1}\ar[r] &
\BZ^{2g-1} \oplus \BZ \ar[r] & \BZ\oplus \BZ }
\end{eqnarray*}

\bigskip

By Poincar\'{e} duality, $H^2(Y;\BZ)$ is also two dimensional;  we
will need an explicit basis for this space as well.  Let
$\Omega_\Sigma$ be the volume form on $\Sigma$.  There exists a
2--form $\Omega$ on $Y$ which restricts to the volume form
$\Omega_\Sigma$ on each fiber $\Sigma$ of $p:Y\ra S^1$. One way to
get this 2--form is to take a metric $g_Y$ on $Y$ such that
$\theta$ is harmonic, then the Hodge dual $\Omega = \star \theta$
is the desired form.

%We could also use the fact that $\varphi^* \Omega_\Sigma =
%\Omega_\Sigma$, extend2--form $\Omega_\Sigma$ can be extended to
%$\Sigma\times [0,1]$ and glued together to get the desired form on
%$Y$.

To get the second basis element, glue up $\alpha_1$ to get a
1--form on $Y$ by setting
$$\alpha(x,t) = \rho(t)\alpha_1(x) + (1-\rho(t))(\varphi^*(\alpha_1(x)))
- (\frac{d}{dt}\rho(t))g(x)dt,$$ where $g\in \Omega^0(Y)$ is the
function given by $\varphi^*(\alpha_1) = \alpha_1 +\beta_1 + dg$.
This 1--form is not closed, in fact
$$d\alpha =  \beta \wedge
(\frac{d}{dt}{\rho}) \theta,$$  but $\alpha\wedge\theta \in
\Omega^2(Y;\BR)$ represents a closed, nontrivial, integral
cohomology class.

\bigskip

The desired 4-manifold $X_g$ is a circle bundle over $Y$ with
Euler class $\chi = [\alpha \wedge \theta]$ in $H^2(Y;\BZ)$. There
is a connection 1--form $\eta$ for this bundle whose curvature
form is $\alpha\wedge \theta$, i.e., $d\eta =
\pi^*(\alpha\wedge\theta)$. Set
$$\omega = \pi^*(\Omega) + \pi^*(\theta)\wedge \eta,$$
then $\omega$ is non-degenerate because $\omega^2$ is the pullback
of the volume form on $Y$ wedge a nowhere-zero 1-form $\eta$. It
is also closed since
$$d\omega = -\pi^*(\theta)\wedge d\eta = 0.$$ Thus $(X_g,\omega)$ is
a symplectic manifold.

\section{Properties of $X_g$}

Next we show that $b_+(X_g)=1$, and that $X_g$ satisfies the
conditions $K\cdot [\omega]>0$ and $K^2=0$.  A calculation using
the Gysin sequence and the fact that $\chi\not=0$ shows that the
cohomology $H^2(X_g;\BZ)$ is 2--dimensional.

\begin{eqnarray*}
\xymatrix{  H^0(Y) \ar[r]^{\cdot \ \cup \chi} \ar@{=}[d] & H^2(Y)
\ar[r]^{\pi^*} \ar@{=}[d] & H^2(X_g) \ar[r] \ar@{=}[d] & H^1(Y)
\ar@{=}[d] \ar[r]^{\cdot \ \cup \chi}
& H^3(Y) \ar[r] \ar@{=}[d] & 0 &\\
 \BZ \ar[r] & \BZ\oplus \BZ \ar[r] & \BZ\oplus\BZ \ar[r]
& \BZ\oplus\BZ \ar[r] & \BZ \ar[r] & 0 }
\end{eqnarray*}

\noindent Hence $\langle [\pi^*(\Omega)], [\pi^*(\theta)\wedge
\eta] \rangle$ is a basis for $H^2(X;\BZ)$ with intersection
matrix
$$Q_X = \left( \begin{array}{cc} 0 & 1\\ 1 & d\end{array} \right),$$
for some $d\in\BZ$. ($[\pi^*(\Omega)]^2=0$ by naturality of
$\pi^*$.) This quadratic form clearly has signature zero, implying
that $b_+(X)=1$.

Next we define a compatible complex structure $J$ and use it to
calculate $K$.  Because $\Sigma$ intersects trivially with the
Poincar\'e dual of $\chi$, the tangent space $T\Sigma$ lifts to a
subspace of $TX$ at each point in $X$.  Fix a complex structure on
the $\Sigma$. On the subspace $T\Sigma\subset TX$, define $J$ to
be the complex structure given by $\Sigma$. Also at each point in
$X$, there are two natural vectors, $\frac{\partial}{\partial t}$
and $\T$, given by $\theta(\frac{\partial}{\partial t}) =1$ and
$\eta(\T)=1$, which are linearly independent to each other and the
subspace $T\Sigma$. Define $J \frac{\partial}{\partial t} = \T$
and $J T = -\frac{\partial}{\partial t}$ on the span of these
vectors.  It is easy to see that $J$ is compatible with $\omega$.
With this complex structure the tangent space of $TX$ splits into
a direct sum and
$$K=-c_1(TX,J) = -c_1(T\Sigma\oplus \underline{\BC}) =
(2g-2)[\pi^*\Omega].$$ (In Section 6 we will see another proof of
this equality.) Hence $K^2=0$ and $K\cdot [\omega] = 2g-2 \geq 0$.

We are interested in $g{>}1$, but the discussion also holds for
the case when $g{=}1$. In that case the construction yields a
symplectic version of a Calabi-Yau 4--manifold, i.e., a symplectic
manifold with $K=0$.

\section{$X_g$ does not admit a complex structure}

\label{sec:complex_structure}

Recall that a closed 4--manifold with even $b_1$ is complex if and
only if it is \Kahler\ (cf. \cite{sw:4man_kirby_calc}).  By the
Gysin seqence, $\pi^*:H^1(Y;\BZ) \ra H^1(X_g;\BZ)$ is an
isomorphism, and therefore $b_1(X_g){=}2$. Thus to show that $X_g$
is not complex it is enough to show that it is not of Lefschetz
type.

Suppose that $X_g$ is \Kahler.  Then for some $L\in H^2(X_g;\BZ)$,
$$\cdot \cup L: H^1(X_g;\BZ) \ra
H^3(X_g;\BZ)$$ is an isomorphism by the Hard Lefschetz Theorem.
This in turn implies that there is a non-degenerate quadratic
form,
$$q:H^1(X_g;\BZ)\ot H^1(X_g;\BZ) \ra \BZ,$$ given by $q(a,b) = a\cup
b\cup L$.

We will show that the class $[\pi^*(\theta)]$  annihilates
$H^1(X_g;\BZ)$ by cup product, contradicting that $q$ is
non-degenerate. Notice that this contradiction also means that
$X_g$ is not of Lefschetz type.

Instead of working with the cup product structure of $X_g$, we can
work with $\theta$ and the cup product structure on $Y$ instead.
This is because of the naturality of the cup product and the
isomorphism $H^1(Y;\BZ)\cong H^1(X_g;\BZ)$ given by the Gysin
sequence.

Apply $[\theta]$ to the basis elements $\{ [\theta], [\beta] \}$
described above.  Clearly  $$[\theta]\cup [\theta] =
[\theta\wedge\theta] =0$$ since $\theta$ is a closed 1--form.  To
show that $[\theta]\cup [\beta] = 0$, evaluate the 2--form
$\theta\wedge \beta$ on a basis of $H_2(Y;\BZ)$.  Consider
$\Sigma$ as one of the fibers of the fibration $p:Y\ra S^1$, then
$$\langle [\theta]\cup [\beta],[\Sigma]\rangle \ = \ \int_\Sigma
\theta \wedge \beta \ = \ 0.$$  The other basis element of
$H_2(Y;\BZ)$ is represented by $a_1\times S^1$ --- the torus that
is transverse to the fiber $\Sigma$ in $Y$. (Recall that $a_1$ is
the loop such that $\varphi_*(a_1)=a_1$.)  Then
$$\langle [\theta]\cup [\beta],[a_1\times S^1]\rangle \ = \
\int_{a_1\times S^1} \theta \wedge \beta \ = \ \beta_1(a_1) \ = \
0.$$

Thus $X_g$ is not of Lefschetz type and does not admit a complex
structure.

\section{Seiberg-Witten invariants and the wall crossing formula}
\label{sec:sw_inv}

The Seiberg-Witten invariants of a 4--manifold $X$ for a fixed
\Spinc\ structure is the count of signed solutions to a partial
differential equation on the spinor bundle. In this section we
give a brief introduction of the relevant topics needed to
understand why the Seiberg-Witten invariants are diffeomorphism
invariants for $X_g$ and to understand why they do not carry
metrics of positive scalar curvature.

Fix a metric $h$ on $X$.  A \Spinc\ structure on the 4-manifold
$X$ is a pair $\xi = (W, \sigma)$ consisting of a rank 4 complex
bundle $W$ with a hermitian metric (the spinor bundle) and an
action $\sigma$ of 1-forms on spinors,
\[\sigma: T^*X \ra \mbox{End}(X),\]
which satisfies the property that, if $e^1, e^2, e^3, e^4$ are an
orthonormal coframe at a point in $X$, then the endomorphisms
$\sigma(e^i)$ are skew-adjoint and satisfy the Clifford relations
\[\sigma(e^i)\sigma(e^j) + \sigma(e^j)\sigma(e^i) = -2\delta_{ij}.\]
The bundle $W$ decomposes into two bundles of rank $2$, $W=W^+
\oplus W^-$, with $\det W^+ = \det W^-$.  The bundle $W^-$ is the
subspace annihilated by the action of self-dual $2$-forms
($\sigma$ can be extended to an action of 2--forms).

By coupling it with a $U(1)$-connection $A$ on $\det W^+$ with the
Levi-Civita connection we can define a connection on $W^+$. Use
this connection to define a Dirac operator
$\dirac^+_A:\Gamma_X(W^+)\ra \Gamma_X(W^-)$ from the space of
smooth sections of $W^+$ to $W^-$. The $4$-dimensional perturbed
Seiberg-Witten equations for a section $\Psi \in \Gamma_X(W^+)$
and a $U(1)$-connection $A$ on $\det W^+$ are:

$$
\label{eq:sw_4dim}
\begin{array}{r@{=}l}
F^+_A + \delta  - q(\Psi) \ \ & \ \ 0,\\
\dirac^+_A(\Psi) \ \ & \ \ 0.
\end{array}
$$
Here $F^+_A$ is the projection of the curvature onto the self-dual
two forms, $\delta \in \Omega^+(X;\I\BR)$ is self-dual $2$-form
used to perturb the equations, and $q:~\Gamma_X(W^+) \ra
\Omega^+(X, \I\BR)$ defined by ${q(\Psi) = \Psi \ot \Psi^* -
\frac12 |\Psi|^2}$ is the adjoint of Clifford multiplication by
self-dual $2$-forms,i.e,
$$
\label{eq:adjointprop}
 \langle \sigma(\beta) \Psi, \Psi\rangle_{W^+} \ = \ 4 \langle \beta, q(\Psi)\rangle_{\I\Lambda^+}
$$
for all self-dual $2$-forms $\beta \in\Omega^+(X;\I\BR)$ and all
sections $\Psi$.

For a fixed metric $h$ and perturbation term $\delta$, the moduli
space $\mathcal{M}(X,\xi, h,\delta)$ is the space of solutions to
the Seiberg-Witten equations modulo the action of the gauge group
$\mathcal{G}=Map(X,S^1)$. For a generic choice of $(h, \delta)$
the moduli space is a compact, oriented, smooth manifold of
dimension
$$
\label{eq:dim_count} d(\xi)= \frac14\left( c_1(\xi)^2 - 2 \chi(X)
-3\sigma(X)\right)
$$
which is independent of metric and perturbation when $b_+(X)>1$.

The Seiberg-Witten invariant $SW_X(\xi)$ is a suitable count of
solutions. Fix a base point in $M$ and let $\mathcal{G}^0 \subset
Map(X,S^1)$ denote the group of maps which map that point to 1.
The based moduli space, denoted by $\mathcal{M}^0$, is the
quotient of the space of solutions by $\mathcal{G}^0$. When the
moduli space $\mathcal{M}(X,\xi,h,\delta)$ is smooth,
$\mathcal{M}^0$ is a principal $S^1$-bundle over
$\mathcal{M}(X,\xi,h,\delta)$. For a given \Spinc\ structure
$\xi$, the 4-dimensional Seiberg-Witten invariant $SW_X(\xi)$ is
defined to be $0$ when $d(\xi)<0$, the sum of signed points when
$d(\xi)=0$, or if $d(\xi)>0$, it is the pairing of the fundamental
class of $\mathcal{M}(X,\xi,h,\delta)$ with the maximal cup
product of the Euler class of the $S^1$-bundle $\mathcal{M}^0$.

Let $G(X)$ be the product space of metrics and
$\Gamma(\Lambda_+)$. A pair  $(h,\delta) \in G(X)$ is called a
{\em good pair} if the moduli space $\mathcal{M}(X,\xi,h,\delta)$
is a smooth manifold without reducible solutions (ie., solutions
$(A,\Psi)$ where $\Psi\equiv 0$). The bad pairs form a `wall'
inside $G(X)$ of dimension $b_+(X)$.   When $b_+>1$ a cobordism
can be constructed between any two moduli spaces of good pairs,
making the Seiberg-Witten invariants independent of metric and
perturbation. However, when $b_+(X){=}1$ it is possible that two
good pairs cannot be connected through a generic smooth path in
$G(X)$ without crossing a wall of bad pairs where reducible
solutions occur. Passing through a bad pair could cause a
singularity to occur in the cobordism. For a general $b_+=1$
manifold, this will often break the smooth invariance of the
Seiberg-Witten invariant.

In 1995 Li and Liu proved a general wall crossing formula which
describes how the  Seiberg-Witten invariants change when crossing
a wall \cite{sw:wallcross}.  The wall crossing formula for $X_g$
can be stated as follows:

\begin{theorem}[Li-Liu] Let $\mathcal{\xi}\in H^2(X_g;\BZ)$ a \Spinc\
structure with $d(\xi) \geq 0$. There exists a basis $\{y_1,y_2\}$
of $H^1(X_g;\BZ)$ depending on $\xi$  such that $$\langle y_i \cup
y_i \cup \frac{\xi}{2}, [X_g]\rangle \ = \ 0$$ for $i=1,2$. Then,
after crossing a wall, the $SW_X(\xi)$ changes by $\pm y_1y_2
\frac{\xi}{2}[X_g]$.
\end{theorem}

If there exists an element in $H^1(X;\BR)$ which annihilates
$H^1(X;\BR)$ by cup product, the Seiberg-Witten invariant of $X$
does not change after crossing any wall in $G(X)$. In that case
the Seiberg-Witten invariants of $X$ are independent of the metric
and perturbation, and are smooth invariants of the manifold.

We finish this section by showing that $X_g$ does not carry a
metric of positive scalar curvature.  Since $X_g$ is symplectic,
the Seiberg-Witten invariants are nontrivial by a theorem of
Taubes \cite{sw:sw_inv_and_symp_form}. It was shown in Section
\ref{sec:complex_structure} that $[\pi^*(\theta)]$ annihilates
$H^1(X_g;\BR)$, thus the Seiberg-Witten invariants for a given
\Spinc\ structure are the same in both chambers by the wall
crossing formula above. The existence of an irreducible solution
for all good pairs then implies the standard fact in
Seiberg-Witten theory that $X_g$ can not have a metric of positive
scalar curvature (cf. \cite{sw:4man_kirby_calc}).

\section{Seiberg-Witten invariants of $X_g$}

In this section we explicitly calculate the  Seiberg-Witten
invariants of the manifold $X_g$ for any genus $g>0$ using the
formula described in \cite{sw:sworbcirc}.

First we calculate the Milnor torsion of the three manifold $Y$
directly from the fundamental group of $Y$,
$$\pi_1(Y) = \langle t, a_1,b_1,  \ldots a_g, b_g \ | \
\prod [a_i,b_i], \ t a_i t^{-1} \varphi_*^{-1}(a_i), \ t b_i
t^{-1} \varphi_*^{-1}(b_i)\rangle. $$ Here the loop $t$ is  $pt
\times [0,1]$ in $Y =\Sigma \times [0,1] / \sim$. The Alexander
matrix, calculated by Fox calculus from the group presentation
above, is equal to:
\[
\left(\begin{array}{ccc} 0 & 1-b_1 & 0 \\
0 & t-1 & 0\\
1-b_1 & -b_1 & t-1 \end{array}\right) \bigoplus \underbrace{K_8
\oplus \cdots \oplus K_8}_{g-1},
\]
where $K_8=\left(\begin{array}{cc} t-1 & -1\\ -1 &
t-2\end{array}\right)$.  The Alexander polynomial is the
polynomial which generates the first elementary ideal of this
matrix.  Symmetrizing this polynomial gives the Milnor torsion of
$Y$,
$$\Delta_Y(t) = (t^{-1} -3 + t)^{g-1}.$$
By a theorem of Meng-Taubes \cite{sw:swequalmilnortorsion}, the
Seiberg-Witten invariant of $Y$ is
$$ \mathcal{SW}_Y(\mathfrak{t}) = \Delta_Y(\mathfrak{t}^2) = (\mathfrak{t}^{-2}-3+\mathfrak{t}^2)^{g-1}$$
where $\mathfrak{t} = \exp(PD[t]) = \exp(\Omega_\Sigma)$.

By the formula in \cite{sw:sworbcirc}, the Seiberg-Witten
invariants of $X_g$ are found by adding the coefficients of
$\mathcal{SW}_Y$ for all \Spinc\ structures which differ by a
multiple of $\chi = [\alpha \wedge \theta]$.  Since \Spinc\
structures in $\mathcal{SW}_Y$ are multiples of $\Omega_\Sigma$,
the polynomial does not change, except that we need to replace
$\mathfrak{t}$ with $\mathfrak{s} = \exp(\pi^*(\Omega_\Sigma))$.

Note that this calculation confirms that $K=
(2g-2)\pi^*[\Omega_\Sigma]$. (Taubes showed that the canonical
class in this case is the top power of $\mathcal{SW}_{X_g}$.) It
also shows that $X_g$ is not diffeomorphic to $X_h$ for $g\not=h$
although this already follows from the calculation of the
fundamental group.

%\section{Conclusion}

%There are other diffeomorphisms than $\varphi$ that one can use to
%construct symplectic 4--manifolds.  Such maps might lead to
%symplectic manifolds with interesting properties.  In particular,

%\begin{quote} \em Does there exists an infinite number of
%surface diffeomorphisms with the property that each diffeomorphism
%$\varphi:\Sigma_g\ra\Sigma_g$ satisfies ${\dim \ker (\varphi^* -
%Id) =1}$ and $Y = \Sigma\times [0,1]/\sim$ has Milnor torsion
%equal to 1?
%\end{quote}

%Such manifolds might then be used to construct new examples of
%symplectic, non-\Kahler\ 4--manifolds which behave like Calabi-Yau
%manifolds.  Right now the only known examples are special
%$T^2$-bundles over $T^2$.

%**************************************************
%******* Bibliography Stuff ***********************
%**************************************************

%**************************************************

%******* End of document **************************

%**************************************************

\end{document}